\documentclass[a4paper,12pt]{amsart}
\usepackage{brunnian}
\usepackage[labelsep=colon]{caption}
\usepackage{float}
\restylefloat{figure}
\usepackage{graphicx}
\usepackage{harvard}
\usepackage{hyperref}
\usepackage{mathrsfs}
\usepackage[hang]{subfigure}
\usetikzlibrary{arrows,%
  decorations.markings,%
  decorations.pathmorphing,%
  matrix,%
  shapes}
\usepackage{wrapfig}

\hypersetup{
  pdftitle={Higher order architecture of collections of objects},%
  pdfauthor={Nils A. Baas},%
}

\DeclareMathOperator{\Bond}{Bond}
\DeclareMathOperator{\Coll}{Coll}
\DeclareMathOperator{\Glob}{Glob}
\DeclareMathOperator{\id}{id}
\DeclareMathOperator{\Mor}{Mor}
\DeclareMathOperator{\Pre}{Pre}
\DeclareMathOperator{\Sh}{Sh}
\DeclareMathOperator{\source}{source}
\DeclareMathOperator{\target}{target}
\newcommand{\A}{\mathcal{A}}
\newcommand{\C}{\mathscr{C}}
\newcommand{\cone}{\mathcal{C}}
\newcommand{\diagram}[3]{\matrix (#1) [matrix of math nodes,row
  sep={#2},column sep={#3},text height=1.5ex,text depth=0.25ex]}

\newcommand{\fusionscript}{\square_\HH}
\newcommand{\G}{\mathscr{G}}
\newcommand{\HH}{\mathscr{H}}
\newcommand{\Mfd}{\mathit{Mfd}}
\newcommand{\PP}{\mathcal{P}}
\newcommand{\PPP}{\mathscr{P}}
\newcommand{\Sets}{\mathit{Sets}}
\newcommand{\Shyper}{\mathscr{S}}
\newcommand{\Sieve}{\mathcal{S}}
\newcommand{\Strat}{\mathit{Strat}}
\newcommand{\T}{\mathcal{T}}

\theoremstyle{definition}
\newtheorem{defn}{\bfseries Definition}
\newtheorem*{definition}{\bfseries Definition}
\newtheorem*{idea}{\bfseries The basic idea}
\newtheorem*{thp}{\bfseries The Hyperstructure Principle}

\makeatletter
\g@addto@macro\th@definition{\thm@headpunct{:}}
\makeatother

\tikzset{>=stealth}

\title{Higher order architecture of collections of objects}
\author[N.A. Baas]{Nils A.\ Baas$^\ast$}
\dedicatory{Deparment of Mathematical Sciences, NTNU, N-7491
  Trondheim, Norway}
\thanks{$^\ast$Email: baas@math.ntnu.no}

\begin{document}

\begin{abstract}
  We show that on an arbitrary collection of objects there is a wide
  variety of higher order architectures governed by hyperstructures.
  Higher order gluing, local to global processes, fusion of
  collections, bridges and higher order types are discussed.  We think
  that these types of architectures may have interesting applications
  in many areas of science.\\

\end{abstract}

\maketitle

\section{Introduction}
\label{sec:introduction}

In all human activites we consider objects and collections of
objects.  We study how we can make new objects and collections out of
old ones by using their properties, relations and interactions.  In
order to do so it is important to have in mind what kind of
architectures to use when forming new structures and organizations out
of given collections of objects.  What are the possibilities?

The main purpose of this paper is to describe new architectures and
discuss how to organize a collection of objects in theory and
practice.  We will introduce higher order architectures which lead to
new and unexplored structures and organizations extending for example
categorical organizations.  We basically follow the hyperstructure
idea with a new extension.  We will discuss examples and situations
where these new architectures may be useful.

Our main purpose is to point out that there is a plethora of higher
order architectures waiting to be explored.

\section{Architectures}
\label{sec:architectures}

When we are given a collection of objects and we want to organize
them, what are the architectures we may use?

Let us think of an arbitrary collection $\C$ as a set of objects.  In
settings where we would like to avoid set theory the general idea will
apply as well.

Finite collections can be organized in many ways.  Given a finite set
$X$ we may consider the elements or objects as vertices in a graph or
a network if we assign weights on the edges.  We may extend this to
hypergraphs where the edges are replaced by subsets of objects instead
of pairs.

Another useful organization of a set $X$ is as a simplicial complex
which is specified by a set of finite non-empty subsets $S \subset
\PP(X)$ (the power set of $X$, i.e.\ the set of all subsets) such that
if $s\in S$ and $t \subset s$ ($t$ non-empty), then $t\in S$ and $x\in
X$ implies $\{x\} \in S$.

A topological space is an organization of a set $X$ which is given by
a subcollection $\T \subset \PP(X)$ satisfying certain axioms.  The
elements of $\T$ are called open sets.

Similarly a measure space is defined by another type of subcollection
$\A \subset \PP(X)$ satisfying a different set of axioms.  Many types
of structures can be described in this way as in
\citeasnoun{Bourbaki}.

A category is an organization on a collection of objects (class or
set) where one assigns a set
\begin{equation*}
  \Mor(X,Y)
\end{equation*}
to every ordered pair of objects.  In addition there is a composition
law.

We will define an organizational architecture which encompasses these
examples.

\begin{idea}
  Instead of assigning a set $\Mor(X,Y)$ to every ordered pair of
  objects we will assign a set of bonds to any collection of objects
  --- finite, infinite or uncountable:
  \begin{equation*}
    \Bond(X,Y,Z,\ldots) \quad \text{or} \quad \Bond(c\in \C)
  \end{equation*}
  $\C$ being a collection or parametrized family of objects.  We may
  also consider ordered collections or collections with other
  additional properties.  Bonds extend morphisms in categories and
  higher bonds create levels and extend higher morphisms (natural
  transformations and homotopies, etc.).  This will be the basis for
  the creation of new global states.\\

  More formally our extended architecture is as follows:

  $X$ is a set representing an arbitrary collection of objects.  Then
  we consider the iterated power sets
  \begin{equation*}
    X, \PP(X), \PP^2(X),\ldots,\PP^k(X),\ldots
  \end{equation*}
  For $S \in \PP^k(X)$ for some $k$ we assign a set of properties or
  states as follows:
  \begin{equation*}
    \Omega \colon \PP^k(X) \to \Sets
  \end{equation*}
  and we think of $\Omega$ as a presheaf functor in a suitable way.
  Let us put
  \begin{equation*}
    \PPP(X) = \bigcup_{k = 0}^\infty \PP^k(X)
  \end{equation*}
  and extend
  \begin{equation*}
    \Omega \colon \PPP(X) \to \Sets.
  \end{equation*}
  The philosophy is that we first pick out the generalized subset we
  want to bind, then we assign the properties that will be involved in
  the binding process ending up with pairs
  \begin{equation*}
    (S,\omega) \quad S\in \PPP(X) \text{ and } \omega \in \Omega(S).
  \end{equation*}
  In general we may also consider situations where $\PPP(X)$ is a
  suitable space of collections of spaces, algebras, etc.  This means
  that $\PPP(X)$ may be a space where for example the points are
  spaces as well, like in moduli type spaces.

  A typical $S$ may look like the configuration in Figure
  \ref{fig:generalized_collection}.
  \begin{figure}[H]
    \centering
    \includegraphics[scale=0.65]{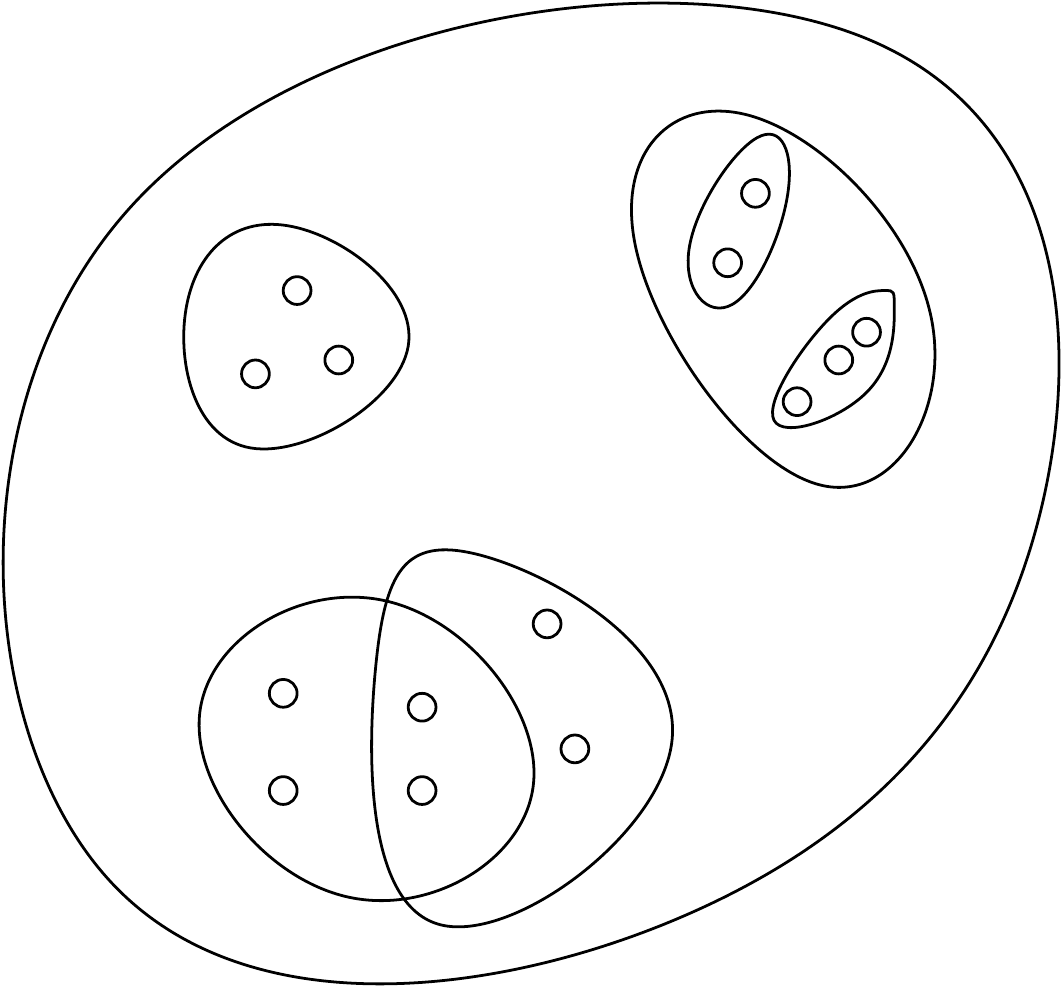}
    \caption{A typical generalized collection}
    \label{fig:generalized_collection}
  \end{figure}
  These pairs are now our building blocks for new collections and we
  have to specify how we bind them together.  Formally as for the
  original hyperstructures we define
  \begin{equation*}
    \Gamma = \{ (S,\omega) \mid S\in \PPP(X) \text{ and } \omega \in
    \Omega(S)\}
  \end{equation*}
  $B$ is then an assignment of a set of bonds
  \begin{equation*}
    B \colon \Gamma \to \Sets
  \end{equation*}
  which we also may think of as a presheaf functor in a suitable way.
  Clearly some binding sets may be empty if the collection does not
  bind under the present circumstances.

  This extends notions of morphisms, relations, etc.
\end{idea}

\section{Higher order architectures}
\label{sec:higher}

What about using the bonds of collections as building blocks for new
bonds of collections of collections?

This is a natural extension in order to create higher order structures
similar to higher power sets --- $\PP^k(X)$, relations of relations,
morphisms of morphisms in higher category theory.

Let us consider an example.  Suppose that we are given a finite set of
objects $V$.  In a given context we want to maximize the number of
interactions in $V$.  The optimal structure is clearly the complete
graph, which may not always be obtainable, so in general we end up
with a subgraph of the complete graph.  Does the evolutionary process
of creating more and more interactions stop here?

No, because we may let subsets of $V$ --- not just pairs --- interact.
When we have exhausted this possibility, we may proceed to subsets of
subsets, etc.  This shows that even on a finite set there is an
enourmous number of possibilities of evolving new interactions.  We
use this example as a guide for the general contruction.

In describing higher order architectures of organizations of
collections we follow the ``General Principle'' of hyperstructures in
\citeasnoun{HAM}, for background material see \cite{EHH,Hnano}.  The
objects may be abstract or physical.

This means that we start with a basic collection $\C$.  As in
\cite{HAM,NSCS,SO} we set $\C = X_0$, and then form
successively:
\begin{equation*}
  X_0, \quad \Omega_0, \quad \Gamma_0, \quad B_0.
\end{equation*}
In order to create the next level we put
\begin{equation*}
  X_1 = \{ b_0 \mid b_0 \in B_0(S_0,\omega_0), S_0 \in \PPP(X_0)
  \text{ and } \omega_0 \in \Omega_0(S_0)\}.
\end{equation*}
Depending on the situation we now can choose $\Omega_1$ and $B_1$
according to what we want to construct or study and then repeat the
construction.

This is not a recursive procedure since new properties and bonds arise
at each level.

Hence a higher order architecture of order $n$ is described by:
\begin{equation*}
  \HH_n \quad \colon \quad \begin{cases}
    X_0, \quad \Omega_0, \quad \Gamma_0, \quad B_0\\
    X_1, \quad \Omega_1, \quad \Gamma_1, \quad B_1\\[0.25cm]
    \hspace*{2cm} \vdots\\[0.25cm]
    X_n, \quad \Omega_n, \quad \Gamma_n, \quad B_n.
  \end{cases}
\end{equation*}
At the technical level we require that
\begin{equation*}
  B_i(S_i,\omega_i) \cap B_i (S_i',\omega_i') = \emptyset
\end{equation*}
for $S_i \neq S_i'$ (``a bond knows what it binds'') in order to
define the $\partial_i$'s below, or we could just require that the
$\partial_i$'s exist.

The level architectures are connected by ``boundary'' maps as follows:
\begin{equation*}
  \partial_i \colon X_{i + 1} \to \PPP(X_i)
\end{equation*}
defined by
\begin{equation*}
  \partial_i(b_{i + 1}) = S_i
\end{equation*}
and maps
\begin{equation*}
  I_i \colon X_i \to X_{i + 1}
\end{equation*}
such that $\partial_i \circ I_i = \id$.  $I_i$ gives a kind of
``identity bond''.

The extensions allowing bindings of subsets or subcollections of
higher power sets add many new types of architectures of
hyperstructures.  See \citeasnoun{HTD} for examples.

The structures $\HH_n$ is an extension of what we called
hyperstructures in \cite{HAM,NSCS,SO} but we cover them by the same
name.  In order to understand higher order structures we introduce,
describe, define and detect hyperstructures.

By selective choices of observational mechanisms of properties and
states the combinatorics may be kept at a reasonable level.

The general higher order architectural pattern is:

Start with a collection $\C = \C_0$,
\begin{equation*}
  \begin{array}{r@{\ }c@{\ }l}
    \C_0 & = & \Coll_0\\
    \C_1 & = & \Bond_0\left(\Coll_0\right)\\[0.25cm]
    & \vdots &\\[0.25cm]
    \C_n & = & \Bond_{n - 1}\left(\Coll_{n -
        1}\left(\ldots\left(\Coll_0\right)\ldots\right)\right)
  \end{array}
\end{equation*}
where $\Coll_k \subset \C_k$ has been selected.

\begin{definition}
  A higher order architecture of order $n$ on a collection of objects
  $\C$ is given by an $\HH_n$-structure on $\C$.
\end{definition}

Furthermore, in the realization of such structures in physics,
chemistry, biology, etc.\ there is ``a lot of room'' in the
nano-dimensions as pointed out by R.\ Feynmann.

The conclusion to be made is that there are universes of new higher
order architectures of collections of objects waiting to be explored
both in mathematics and other sciences.

\section{Examples}
\label{sec:examples}

\begin{itemize}
\item[\textbf{I.}] \textbf{Links}\\

  \noindent A link is a disjoint union of embedded circles (or rings)
  in three dimensional space.  They may be linked in many ways.  We
  consider linking as a kind of geometrical or topological binding,
  see Figure \ref{fig:links}.
  \begin{figure}[H]
    \centering
    \subfigure[Brunnian rings]{
      \begin{tikzpicture}[every path/.style={knot},scale=0.5]
        \pgfmathsetmacro{\brscale}{1.8}
        \foreach \brk in {1,2,3} {
          \begin{scope}[rotate=\brk * 120 - 120]
            \colorlet{chain}{ring\brk}
            \brunnianlink{\brscale}{120}
          \end{scope}
        }
      \end{tikzpicture}
      \label{fig:1stBrunnian}
    } \subfigure[2nd order Brunnian rings]{
      \includegraphics[width=0.25\linewidth]{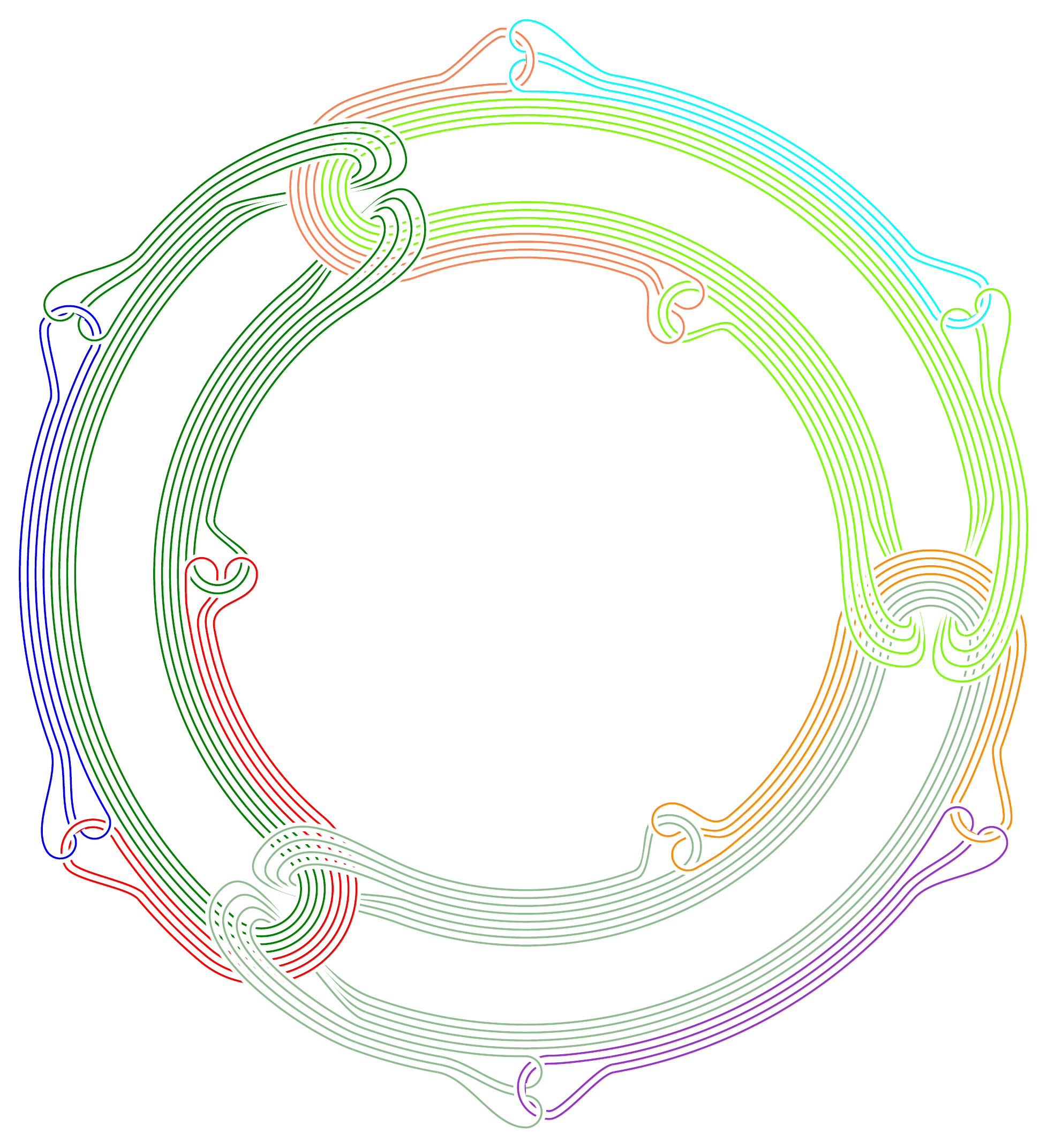}
      \label{fig:2ndBrunnian}
    } \qquad \subfigure[3rd order Brunnian rings]{
      \includegraphics[width=0.25\linewidth]{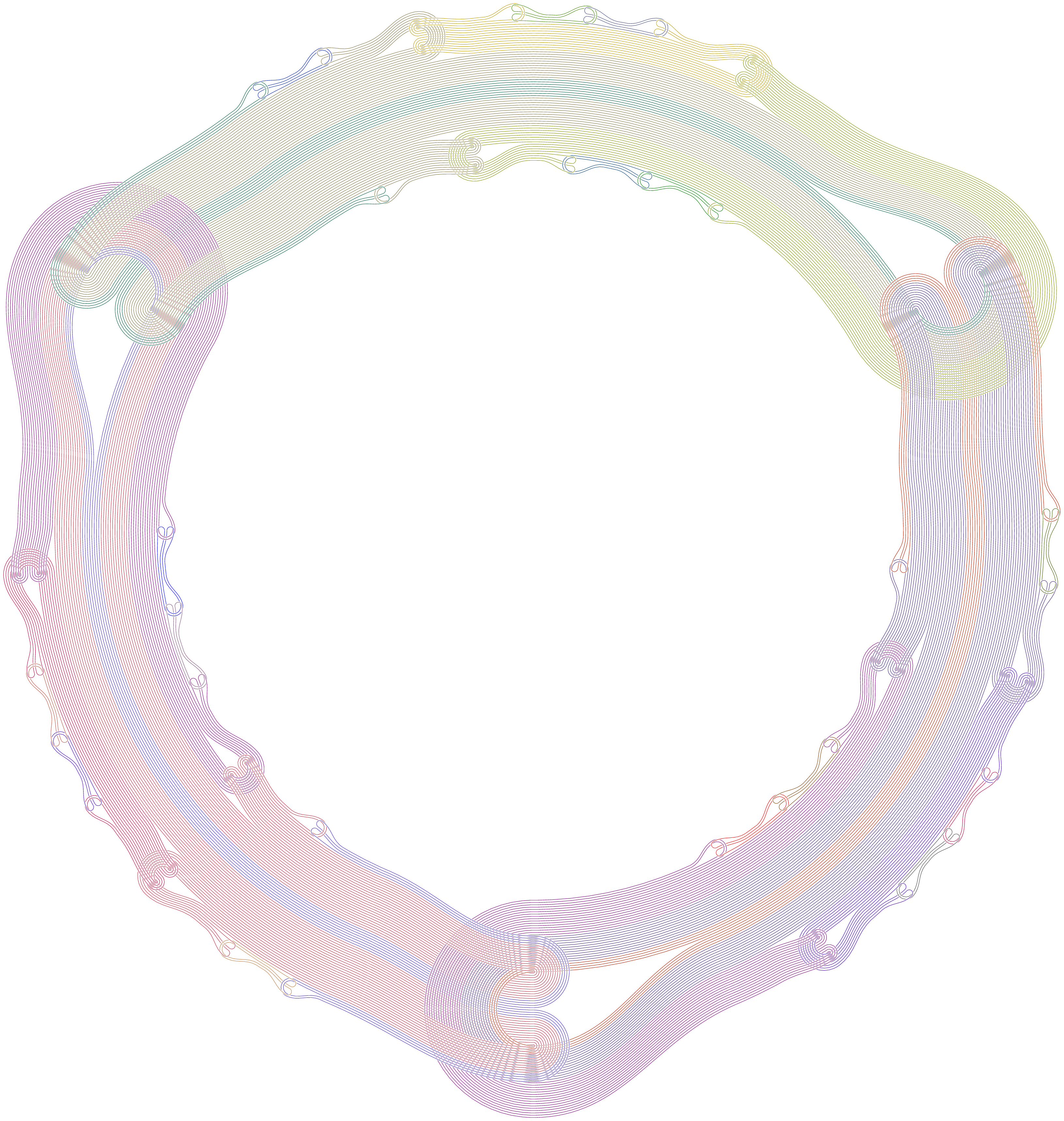}
      \label{fig:3rdBrunnian}
    }
    \caption{Links}
    \label{fig:links}
  \end{figure}
  Once a link has been constructed we may inspect its connectivity
  properties and use these to form links of links --- second order
  links.

  We can then follow the general scheme that we have outlined for
  higher order architecture.  Some of these have been extensively
  studied in \citeasnoun{NS}, especially with respect to Brunnian type
  binding properties.\\
\item[\textbf{II.}] \textbf{Molecules and materials}\\

  \noindent 
  The higher order architectures for geometric links may also act as
  architectures for molecules and how they form materials.  In
  synthetic chemistry a lot of work has been done on the synthesis of
  Borromean rings.  N.\ Seeman was the first to succeed using DNA
  molecules.  In his laboratory work has now been started on the
  synthesis of second order Brunnian rings.  The architectures we have
  outlined may be an interesting and useful guide for future
  synthesis.  A challenging question is what kind of properties such
  materials will have and how they depend on the architecture.  For a
  further discussion of these issues, see \citeasnoun{NSCS}, and
  \citeasnoun{BS}.\\
\item[\textbf{III.}] \textbf{Quantum states}\\

  \noindent Many-body systems play an important role in physics, for
  example in quantum mechanics.  For many-particle systems one may
  introduce the higher order architectures that we have discussed.
  However, can such an architecture be realized as a (bound) quantum
  state?

  \sloppy 
  In 1970, the Russian physicist V.\ Efimov predicted new
  counter-intuitive quantum states where three particles are bound,
  but not two by two.  This is analogous to the Borromean and Brunnian
  linking property.  Experimentally such states were not observed
  until 2006 in ultracold caesium gases.

  This raises the following interesting questions:  Does there exist a
  family of quantum states analogous to higher order Brunnian links?

  May all higher order architectures --- for example for links --- be
  realized as entangled quantum states?  For a further discussion, see
  \citeasnoun{NS} and \citeasnoun{BF}.\\
\item[\textbf{IV.}] \textbf{Geometric bonds}\\

  \noindent The idea is to extend higher cobordisms (cobordisms of
  cobordisms $\ldots$) and nested families of spaces to a
  hyperstructure context.  We refer to \citeasnoun{SO} for a
  discussion of this topic.
\end{itemize}

\section{Higher categories and hypercategories}
\label{sec:hypercat}

In category theory one organizes a collection of objects by assigning
sets of morphisms to ordered pairs.  In our terminology morphisms play
the role of bonds --- binding two objects together.  Then for good
mathematical reasons one introduces morphisms of morphisms ---
two-morphisms, up to general $n$-morphisms, etc.  This is clearly a
special case of the bonds of bonds $\ldots$ situation.

However, it has turned out to be difficult to find a generally
accepted definition of an $n$-category.  From our point of view there
is a recent approach (unpublished) by D.\ Ayala where higher
categories are basically considered as sheaves on manifolds
\begin{equation*}
  \Sh(\Mfd).
\end{equation*}
In this approach the combinatorics of morphisms takes place on
manifolds and the morphisms as bonds may geometrically be thought of
as the cone on two points (source and target).

In our approach we want bonds to bind more than two objects: three,
four,$\ldots$ even whole manifolds.

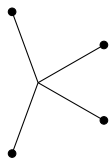
\begin{figure}[H]
  \centering
  \begin{tikzpicture}
    \foreach \a in {30,110,250,330}{
      \filldraw[fill=black] (\a:1cm) circle(0.05cm);
      \draw (\a:1cm) -- (0,0);
    };
  \end{tikzpicture}
  \caption{$\cone(\mathbb{Z}_4) =$ cone on four points}
  \label{fig:singularity}
\end{figure}

The cone in general introduces singularities (Figure
\ref{fig:singularity}).

Intuitively it is a good and natural representation of a bond, and in
general a cone on $P$ is a good example of a bond:
\begin{equation*}
  \cone(P) = \Bond\{p \mid p \in P\}.
\end{equation*}

Let us therefore allow cone-type singularities.  These have been
studied in cobordism theory, \citeasnoun{Sing}, and we may even just
consider arbitrary stratified sets.  They have singular points, but
can be decomposed in a regular way into smooth non-singular pieces
called strata.

Locally any singular point look like the cone on a cone on a cone
$\ldots$ This geometrically iterated cone picture corresponds exactly
to the bonds of bonds of $\ldots$ picture.

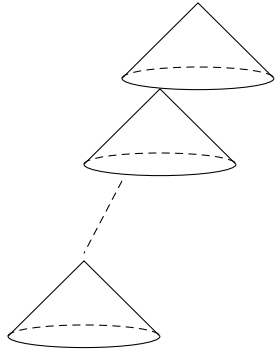
\begin{figure}[H]
  \centering
  \begin{tikzpicture}
    \draw (0,0) -- (1,1) -- (2,0);
    \draw (0,0) arc (180:360:1cm and 0.15cm);
    \draw[densely dashed] (0,0) arc(180:0:1cm and 0.15cm);

    \begin{scope}[xshift=-0.5cm,yshift=-1.14cm]
      \draw (0,0) -- (1,1) -- (2,0);
      \draw (0,0) arc (180:360:1cm and 0.15cm);
      \draw[densely dashed] (0,0) arc(180:0:1cm and 0.15cm);
    \end{scope}

    \begin{scope}[xshift=-1.5cm,yshift=-3.42cm]
      \draw (0,0) -- (1,1) -- (2,0);
      \draw (0,0) arc (180:360:1cm and 0.15cm);
      \draw[densely dashed] (0,0) arc(180:0:1cm and 0.15cm);
    \end{scope}

    \draw[densely dashed] (0,-1.36) -- (-0.5,-2.32);
  \end{tikzpicture}
  \caption{Several levels of cone structure}
  \label{fig:cone_cone}
\end{figure}

This motivates the definition of an extended higher category which we
call a \emph{hypercategory}, basically considered as sheaves on
stratified sets:
\begin{equation*}
  \Sh(\Strat).
\end{equation*}
The details will be developed with D.\ Ayala in future work.

This is an interesting implementation of the more general
hyperstructure architectures.  Cone structures lead to decompositions
and stratifications which again lead to hyperstructures and bonds, see
\citeasnoun{SO}.  Their representations lead to extended field theories.

\section{Compositions of bonds}
\label{sec:compositions}

In the study of collections of objects we emphasize the general notion
of bonds including relations, functions and morphisms.  We get richer
structures when we have composition rules of various types of bonds.
Such compositions should take into account the higher order
architecture giving bonds a level structure.

We experience this situation in higher categories where we want to
compose morphisms of any order.  Suppose that we are given two
$n$-morphisms $f$ and $g$.  They may not be compatible at level $n$
for composition in the sense that
\begin{equation*}
  \target(f) = \source(g).
\end{equation*}
But in a precise way we can iterate source and target maps to get down
to lower levels, and it may then happen that at level $p$ we have
\begin{equation*}
  \target_p^n(f) = \source_p^n(g).
\end{equation*}
Hence composition makes sense at level $p$ and we write the
composition rule as
\begin{equation*}
  \square_p^n
\end{equation*}
and the composed object as
\begin{equation*}
  f \, \square_p^n \, g.
\end{equation*}

In a similar way we can introduce composition rules for bonds in a
general hyperstructure $\HH$.  Let $a_n$ and $b_n$ be bonds at level
$n$ in $\HH$.  Then we get to the lower levels via the boundary maps
\begin{equation*}
  \partial_i \colon X_{i + 1} \to \PPP(X_i)
\end{equation*}
and search for compatibility in the sense that
\begin{equation*}
  \partial_p \circ \cdots \circ \partial_{n - 1} (a_n) = \partial_p
  \circ \cdots \circ \partial_{n - 1} (b_n)
\end{equation*}
or we may just require a weaker condition like
\begin{equation*}
  \partial_p \circ \cdots \circ \partial_{n - 1} (a_n)
  \cap \partial_p \circ \cdots \circ \partial_{n - 1} (b_n) \neq
  \emptyset
\end{equation*}
in order to have a composition defined:
\begin{equation*}
  a_n \, \square_p^n \, b_n
\end{equation*}

For bonds in a hyperstructure we may even compose bonds at different
levels: $a_m$, $b_n$ compatible at level $p$ via boundary maps, allow
us to define
\begin{equation*}
  a_m \, ^m\underset{p}{\square}^n \, b_n
\end{equation*}
as an $m$-bond for $m \geq n$.  Compositional rules are needed and
will appear elsewhere.

Composition may be thought of as a kind of geometric gluing.  We
consider the bonds as spaces, binding collections of families of
subspaces, these again being bonds, etc.  By the ``boundary'' maps we
go down to a level where these are compatible, gluable bond spaces
along which we may glue the bonds within the type of spaces we
consider.

Therefore hyperstructures offer the framework for a new kind of higher
order gluing in which the level architecture plays a major role.  We
will pursue this in the next section.

\section{Higher order gluing and extending sheaves}
\label{sec:gluing}

Hyperstructures are useful tools in passing from local situations to
global ones in collection of objects.  In this process the level
structure is important.  We will here elaborate the discussion of
multilevel state systems in \citeasnoun{SO}.

In mathematics we often consider situations locally at open sets
covering a space and then glue together basically in one stroke ---
meaning there are just two levels local and global, no intermediate
levels.  In many situations dominated by a hyperstructure this is not
sufficient.  We need a more general hyperstructured way of passing
from local to global in general collections.

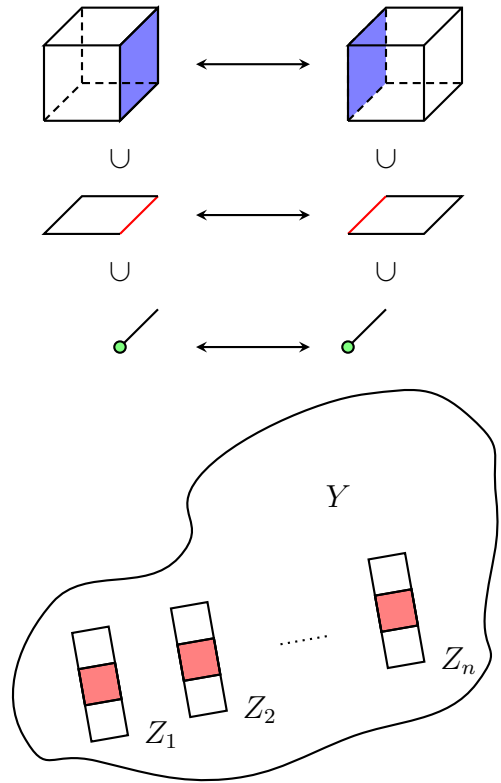
\begin{figure}[H]
  \centering
  \subfigure{
    \begin{tikzpicture}[thick]
      \begin{scope}
        \filldraw[fill=blue!50] (1,0) -- (1.5,0.5) -- (1.5,1.5) -- (1,1)
          -- cycle;
        \draw (0,0) -- (1,0) -- (1.5,0.5) -- (1.5,1.5) -- (0.5,1.5) --
          (0,1) -- (1,1) -- (1.5,1.5);
        \draw (0,0) -- (0,1);
        \draw (1,0) -- (1,1);
        \draw[densely dashed] (0,0) -- (0.5,0.5) -- (0.5,1.5);
        \draw[densely dashed] (0.5,0.5) -- (1.5,0.5);

        \draw[<->] (2,0.75) -- (3.5,0.75);

        \filldraw[fill=blue!50,densely dashed] (4,0) -- (4.5,0.5) --
          (4.5,1.5) -- (4,1) -- cycle;
        \draw (4,0) -- (5,0) -- (5.5,0.5) -- (5.5,1.5) -- (4.5,1.5) --
          (4,1) -- (5,1) -- (5.5,1.5);
        \draw (4,0) -- (4,1);
        \draw (5,1) -- (5,0);
        \draw[densely dashed] (4.5,0.5) -- (5.5,0.5);
      \end{scope}
  
      \begin{scope}[yshift=-0.5cm]
        \node[rotate=90] at (1,0){$\subset$};
        \node[rotate=90] at (4.5,0){$\subset$};
      \end{scope}
      
      \begin{scope}[yshift=-1.5cm]
        \draw[red] (1,0) -- (1.5,0.5);
        \draw (1,0) -- (0,0) -- (0.5,0.5) -- (1.5,0.5);

        \draw[<->] (2,0.25) -- (3.5,0.25);

        \draw[red] (4,0) -- (4.5,0.5);
        \draw (4,0) -- (5,0) -- (5.5,0.5) -- (4.5,0.5);
      \end{scope}

      \begin{scope}[yshift=-2cm]
        \node[rotate=90] at (1,0){$\subset$};
        \node[rotate=90] at (4.5,0){$\subset$};
      \end{scope}
      
      \begin{scope}[xshift=1cm,yshift=-3cm]
        \draw (0,0) -- (0.5,0.5);
        \filldraw[fill=green!50] (0,0) circle(0.075cm);

        \draw[<->] (1,0) -- (2.5,0);

        \draw (3,0) -- (3.5,0.5);
        \filldraw[fill=green!50] (3,0) circle(0.075cm);
      \end{scope}
    \end{tikzpicture}
  } \subfigure{
    \begin{tikzpicture}[scale=0.65,thick]
      \draw[rotate=15] (0.24396864,-1.4967076) .. controls
        (0.49101374,-2.4657116) and (0.665201,-1.8994577)
        .. (1.5439687,-2.3767076) .. controls (2.4227362,-2.8539574) and
        (3.017664,-3.2517743) .. (4.0039687,-3.4167075) .. controls
        (4.9902735,-3.581641) and (5.629609,-3.3427677)
        .. (6.6239686,-3.2367077) .. controls (7.6183286,-3.1306474) and
        (8.599862,-3.4573336) .. (9.383968,-2.8367076) .. controls
        (10.168076,-2.2160816) and (9.585438,-2.3502257)
        .. (9.943969,-1.4167076) .. controls (10.302499,-0.48318952) and
        (10.588963,-0.26462105) .. (10.743969,0.7232924) .. controls
        (10.898975,1.7112058) and (11.202379,1.6006097)
        .. (10.663969,2.4432924) .. controls (10.125558,3.2859752) and
        (9.734352,3.3049438) .. (8.743969,3.4432924) .. controls
        (7.753585,3.581641) and (5.849801,3.2754364)
        .. (5.0439687,2.6832924) .. controls (4.2381363,2.0911484) and
        (5.064207,1.7455099) .. (4.2239685,1.2032924) .. controls
        (3.3837304,0.661075) and (2.3561227,1.5378368)
        .. (1.5239687,0.9832924) .. controls (0.69181454,0.428748) and
        (-0.003076456,-0.5277036) .. (0.24396864,-1.4967076);

      \begin{scope}[xshift=2cm,yshift=-1.5cm,scale=0.75,rotate=10]
        \filldraw[fill=red!50] (0,1) -- (1,1) -- (1,2) -- (0,2) -- cycle;

        \draw (0,0) rectangle(1,3);
        \draw (0,1) rectangle(1,2);
        
        \node[right] at (1.125,0){$Z_1$};
      \end{scope}
      
      \begin{scope}[xshift=4cm,yshift=-1cm,scale=0.75,rotate=10]
        \filldraw[fill=red!50] (0,1) -- (1,1) -- (1,2) -- (0,2) -- cycle;

        \draw (0,0) rectangle(1,3);
        \draw (0,1) rectangle(1,2);
          
        \node[right] at (1.125,0){$Z_2$};
      \end{scope}
        
      \begin{scope}[xshift=6cm,yshift=-0.5cm,rotate=10]
        \draw[dotted] (0,1) -- (1,1);
      \end{scope}
      
      \begin{scope}[xshift=8cm,yshift=0cm,scale=0.75,rotate=10]
        \filldraw[fill=red!50] (0,1) -- (1,1) -- (1,2) -- (0,2) --  cycle;

        \draw (0,0) rectangle(1,3);
        \draw (0,1) rectangle(1,2);
        
        \node[right] at (1.125,0){$Z_n$};
      \end{scope}
        
      \node at (7,3.5){$Y$};
    \end{tikzpicture}
  }
  \caption{Gluing possibility at various levels}
  \label{fig:levels}
\end{figure}

Let us offer two of our intuitions regarding this process.
Geometrically we think of a multilevel nested family of spaces, like
manifolds with singularities represented by manifolds with multinested
boundaries or just like higher dimensional cubes with iterated
boundary structure (corners, edges,$\ldots$).  With two such
structures we may then glue at the various levels of the nesting
(Figure \ref{fig:levels}).
Furthermore, study how states and properties may be ``globalized'',
meaning putting local states coherently together to global states.

Biological systems are put together by multilevel structures from
cells into tissues, organs etc.\ constituting an organism.  Much of
biology is about understanding how cell-states determine organismic
states.  The hyperstructure concept is in fact inspired by biological
systems.

In order to extend the discussion of multilevel state systems in
\citeasnoun{SO} we need to generalize and formulate in a
hyperstructure context the following mathematical notions (see for
example \citeasnoun{MM}):
\begin{itemize}
\item Sieve
\item Grothendieck Topology
\item Site
\item Presheaf
\item Sheaf
\item Descent
\item Stack
\item Sheaf cohomology
\end{itemize}

Let us start with a given hyperstructure
\begin{align*}
  \HH \quad \colon \quad & \{X_0,\ldots,X_{n + 1}\}\\
  & \{\Omega_0,\ldots,\Omega_n\}\\
  & \{B_0,\ldots,B_n\}\\
  & \{\partial_0,\ldots,\partial_n\}\\
\end{align*}
We will now suggest a series of new definitions.

\begin{defn}
  \label{def:sieve}

  A \emph{sieve} on $\HH$ is given as follows: At the lowest level
  $X_0$ a sieve $\Sieve$ on a bond $b_0 (=b_0(S_0,\omega_0))$ is given
  by families of bonds $\{ b_o^{j_0} \}$ (covering families) and
  $b_1$'s such that
  \begin{equation*}
    b_1(\{ b_o^{j_0} \}, b_0),
  \end{equation*}
  $b_0$ may also be replaced by a family of bonds.

  Bond composition with $\{b_0^{j_0}\}$ will produce new families in
  the sieve.


  A \emph{sieve on $\HH$} is then a family of such sieves
  $(\Sieve_k)_{k = 1,\ldots,n}$ --- one for each level.
\end{defn}

We postpone connecting the levels until the definition of a
Grothendieck topology, but this could also have been added to the
sieve definition.

\begin{defn}
  \label{def:Grothendieck_topology}

  A \emph{Grothendieck topology} on $\HH$ is given as follows: First
  we define a Grothendieck topology for each level of bonds.  Consider
  level $0$: to every bond $b_0$ we assign a collection of sieves
  $J(b_0)$ such that
  \begin{itemize}
    \item[i)] (maximality), the maximal sieve on $b_0$ is in $J(b_0)$
    \item[ii)] (stability), let $S \in J(b_0)$, $b_1(b_0',b_0)$, then in
      obvious notation
      \begin{equation*}
        b_1^\ast (S) \in J(b_0')
      \end{equation*}
    \item[iii)] (transitivity), let $S \in J(b_0)$ and $R$ any sieve
      on $b_0,b_0'$ an element of a covering family in $S$, $b_1^\ast
      (R) \in J (b_0')$ for all $b_1$ with $b_1(b_0',b_0)$, then $R
      \in J(b_0)$.
  \end{itemize}
  We call $J(b_0)$ is a $J$-covering of $b_0$.

  This gives a Grothendieck topology for all levels of bonds, and we
  connect them to a structure on all of $\HH$ by defining in addition
  an assignment $J$ of $(b_0,\ldots,b_n)$ where $b_i \in \partial_i
  b_{i + 1}$.

  $J(b_0,\ldots,b_n)$ consists of families of sieves $\{b_o^{j_0}\}
  \in J(b_0),\ldots,\{b_n^{j_n}\}\in J(b_n)$ and bonds
  \begin{equation*}
    \beta_1,\ldots,\beta_{n + 1}
  \end{equation*}
  such that
  \begin{equation*}
    \beta_1(b_0,\{b_o^{j_0}\}),\ldots,\beta_{n +
      1}(b_n,\{b_n^{j_n}\})
  \end{equation*}
  and $b_i^{j_i} \in \partial_i b_{i + 1}^{j_{i + 1}}$.  In a diagram
  we have
  \begin{center}
    \begin{tikzpicture}[descr/.style={fill=white,inner sep=2.5pt}]
      \diagram{d}{2em}{2.5em}{
        b_n & b_{n - 1} & \cdots & b_0\\
        \{b_n^{j_n}\} & \{b_{n - 1}^{j_{n - 1}}\} & \cdots &
        \{b_0^{j_0}\}\\
        J(b_n) & J(b_{n - 1}) && J(b_0).\\
      };

      \path[->,midway,font=\scriptsize]
        (d-1-1) edge node[above]{$\partial$} (d-1-2)
                     edge[<->] node[right]{$\beta_{n + 1}$} (d-2-1)
        (d-1-2) edge node[above]{$\partial$} (d-1-3)
                     edge[<->] node[right]{$\beta_n$} (d-2-2)
        (d-1-3) edge node[above]{$\partial$} (d-1-4)
        (d-1-4) edge[<->] node[right]{$\beta_1$} (d-2-4)
        (d-2-1) edge node[above]{$\partial$} (d-2-2)
                     edge[-,white]
                     node[descr,sloped,text=black]{$\in$} (d-3-1)
        (d-2-2) edge node[above]{$\partial$} (d-2-3)
                     edge[-,white]
                     node[descr,sloped,text=black]{$\in$} (d-3-2)
        (d-2-3) edge node[above]{$\partial$} (d-2-4)
        (d-2-4) edge[-,white]
                     node[descr,sloped,text=black]{$\in$} (d-3-4);
    \end{tikzpicture}
  \end{center}
  Clearly there are many possible choices of Grothendieck topologies,
  and they will be useful in the gluing process and the creation of
  global states.  Examples will be discussed elsewhere, our main point
  here is to outline the general ideas.
\end{defn}

\begin{defn}
  \label{def:site}

  $(\HH,J)$ is called a \emph{hyperstructure site} when
  $J$ is a Grothendieck topology on the hyperstructure $\HH$.
\end{defn}

The categorical notion of a presheaf we modify as follows.  We think
of the $\Omega_i$'s in $\HH$ as state (property, field, dataset, etc.)
\emph{assignments} in a structure preserving way, covering cellular
type (pre)-sheaves as well.  These assignments would be different from
the ones used in the construction of $\HH$.  As in \citeasnoun{SO} we
may let them take values in hyperstructures of states
\begin{equation*}
  \Shyper = \{ \Shyper_0, \Shyper_1,\ldots,\Shyper_n\}
\end{equation*}
$\Shyper_i$ being a hyperstructure such that
\begin{align*}
  \Omega_0 & \text{ takes values in $\Shyper_n$}\\
  \vdots \;\,&\\
  \Omega_i & \text{ takes values in $\Shyper_{n - i}$}\\
  \vdots \;\,&\\
  \Omega_n & \text{ takes values in $\Shyper_0$}
\end{align*}
and we assume that we have bond compatibility of the $\Omega_i$'s and
level connecting assignments $\delta_i$ (``dual'' to the
$\partial_i$'s and acting on collections of bond ``states'') depending
on the Grothendieck topology $J$:
\begin{center}
  \begin{tikzpicture}
    \diagram{d}{2.5em}{2.5em}{
      \Shyper_0 & \Shyper_1 & \cdots & \Shyper_n\\
    };

    \path[->,midway,above,font=\scriptsize]
      (d-1-2) edge[decorate, decoration={snake, amplitude=.4mm,
        segment length=3mm, post length=1mm}] node{$\delta_1$}
        (d-1-1)
      (d-1-3) edge[decorate, decoration={snake, amplitude=.4mm,
        segment length=3mm, post length=1mm}] node{$\delta_2$}
        (d-1-2)
      (d-1-4) edge[decorate, decoration={snake, amplitude=.4mm,
        segment length=3mm, post length=1mm}] node{$\delta_n$}
        (d-1-3);
  \end{tikzpicture}
\end{center}
The $\delta_i$'s may be functional assignments or relational, and the
$\Shyper_i$'s often have an algebraic structure.  In the simplest case
all the $\Shyper_i$'s could just be $\Sets$.

It should be pointed out that the $\Omega_i$'s are already part of the
$\HH$-structure which we now combine with the Grothendieck topology.
This will be useful in constructing new $\Omega$-states of global
bonds.  We consider the $\Omega_i$'s as a kind of ``\emph{level
  presheaves}'' and the $\delta_i$'s giving a kind of ``\emph{global
  matching families}'' --- between levels in addition to levelwise
matching.  However, if we have ``functional'' assignment connectors
$\hat{\delta}_i$'s on $\HH$:
\begin{center}
  \begin{tikzpicture}
    \diagram{d}{2.5em}{2.5em}{
      \Shyper_0 & \Shyper_1 & \cdots & \Shyper_n\\
    };

    \path[->,midway,above,font=\scriptsize]
      (d-1-2) edge node{$\hat{\delta}_1$} (d-1-1)
      (d-1-3) edge node{$\hat{\delta}_2$} (d-1-2)
      (d-1-4) edge node{$\hat{\delta}_n$} (d-1-3);
  \end{tikzpicture}
\end{center}
means that we get a unique state of global bond objects (like an
amalgamation for presheaves but here across levels in addition to
levelwise amalgamation).  Global bonds are ``covered'' as follows (see
\citeasnoun{SO})
\begin{center}
  \begin{tikzpicture}
    \diagram{d}{2.5em}{2.5em}{
      \{ b(i_n) \} & \{ b(i_{n - 1},i_n) \} & \cdots & \{
      b(i_0,\ldots,i_n) \}\\
    };

    \path[->,midway,above,font=\scriptsize]
      (d-1-1) edge node{$\partial_{n -1}$} (d-1-2)
      (d-1-2) edge node{$\partial_{n - 2}$} (d-1-3)
      (d-1-3) edge node{$\partial_0$} (d-1-4);
  \end{tikzpicture}
\end{center}
and states are being levelwise globalized by
\begin{center}
  \begin{tikzpicture}
    \diagram{d}{2.5em}{2em}{
      \Omega_n(\{ b(i_n) \}) & \Omega_{n - 1}(\{ b(i_{n - 1}, i_n)
      \}) & \cdots & \Omega_0 (\{ b(i_0,\ldots,i_n \}).\\
    };

    \path[->,midway,above,font=\scriptsize]
      (d-1-2) edge node{$\hat{\delta}_1$} (d-1-1)
      (d-1-3) edge node{$\hat{\delta}_2$} (d-1-2)
      (d-1-4) edge node{$\hat{\delta}_n$} (d-1-3);
  \end{tikzpicture}
\end{center}

With a slight abuse of notation we write this as
\begin{equation*}
  \Omega \colon (\HH,J) \to \Shyper
\end{equation*}
and define $\Omega = \{\Omega_i\}$ as a \emph{``presheaf''} on
$(\HH,J)$ ($\Pre(\HH,J)$) and when
\begin{equation*}
  \Delta = \{\hat{\delta}_i\}
\end{equation*}
exists we have a unique global bond state.  This is like a
sheafification condition and we call $(\Delta,\Omega)$ a
\emph{globalizer} of the site $(\HH,J)$ with respect to $\Omega$.

$\Omega$ with $\Delta$ extends the sheaf notion here, gluing
within levels and between levels.

\emph{A globalizer is a kind of higher order or hyperstructured sheaf
  covering all the levels.}

The existence of $\Delta$ contains the global gluing data and hence
corresponds to what is often called \emph{descent conditions} and the
hyperstructure collection $\Shyper$ extends the notion of a
\emph{stack} over $\HH$.

Since $\Omega$ is already built into the structure of $\HH$ we may
consider $\Omega$ and $(\Omega, \Delta)$ as ``internal'' presheaves
and sheaves.  However, we may alternatively define an ``external''
assignment
\begin{equation*}
  \Lambda \colon (\HH,J) \to \Shyper
\end{equation*}
where $\Lambda = \{\Lambda_i\}$ and
\begin{equation*}
  \begin{array}{c@{\ }r@{\ }c}
    \Lambda_0 & \text{takes values in} & \Shyper_n\\
    \vdots && \vdots\\
    \Lambda_i & \text{takes values in} & \Shyper_{n-1}\\
    \vdots && \vdots\\
    \Lambda_n & \text{takes values in} & \Shyper_0.
  \end{array}
\end{equation*}
$\Lambda$ is being globalized in the same way as $\Omega$, the
``boundary decomposition''
\begin{center}
  \begin{tikzpicture}
    \diagram{d}{2.5em}{2em}{
      \{ b(i_n) \} & \{ b(i_{n - 1}, i_n)\} & \cdots & \{
      b(i_0,\ldots,i_n)\}\\
    };

    \path[->,midway,above,font=\scriptsize]
      (d-1-1) edge node{$\partial_{n - 1}$} (d-1-2)
      (d-1-2) edge node{$\partial_{n - 2}$} (d-1-3)
      (d-1-3) edge node{$\partial_0$} (d-1-4);
  \end{tikzpicture}
\end{center}
gives rise to
\begin{center}
  \begin{tikzpicture}
    \diagram{d}{2.5em}{2em}{
      \Lambda_n(\{ b(i_n) \}) & \Lambda_{n - 1}(\{ b(i_{n - 1}, i_n)
      \}) & \cdots & \Lambda_0 (\{ b(i_0,\ldots,i_n \}).\\
    };

    \path[->,midway,above,font=\scriptsize]
      (d-1-2) edge node{$\delta'_1$} (d-1-1)
      (d-1-3) edge node{$\delta'_2$} (d-1-2)
      (d-1-4) edge node{$\delta'_n$} (d-1-3);
  \end{tikzpicture}
\end{center}
When
\begin{equation*}
  \Lambda \colon (\HH,J) \to \Shyper
\end{equation*}
is an ``external'' presheaf and $\Delta = \{\delta_i'\}$ exists we
define $(\Delta, \Lambda)$ to be an ``external'' globalizer.

To a site $(\HH,J)$ we may also form a hyperstructure in a suitable way
\begin{equation*}
  \Glob(\HH,J)
\end{equation*}
of globalizers (internal or external) corresponding to a Grothendieck
topos.  If $\HH$ is just a one level categorical structure and
$\Shyper = \Sets$, then
\begin{equation*}
  \Glob(\HH,J) = \Sh(\HH,J).
\end{equation*}
$\HH$ being an $n$-category of some kind is another special case.

This framework is very useful in formulating and studying \emph{local
  to global situations}, and extends the discussion started in
\citeasnoun{SO}, where we as examples studied globalizers given by
compositions of mappings.

In view of our discussion in section \ref{sec:hypercat} this leads to
the following chain of structures
\begin{center}
  \begin{tikzpicture}
    \diagram{d}{2.5em}{2.5em}{
      \Sh(\Mfd) & \Sh(\Strat) & \Glob(\HH\text{-site}).\\
    };

    \path[right hook->]
      (d-1-1) edge (d-1-2)
      (d-1-2) edge (d-1-3);
  \end{tikzpicture}
\end{center}

The concepts introduced here also paves the way for new types of
cohomology of the basic collection or set, say $Z$, with
``coefficients'' in ``presheaves'' and ``globalizers'':
\begin{equation*}
  H^\ast (Z;\PPP) \quad \text{and} \quad H^\ast (Z;\G)
\end{equation*}
$\PPP\in \Pre(\HH,J)$, $\G\in \Glob(\HH,J)$ --- all to be defined in a
suitable setting with $\Shyper$ having an algebraic structure, to be
discussed elsewhere.

In a hyperstructure the ``boundary'' maps $\partial_i$ connect the
levels.  A hyperstructure on a set or collection $Z$, $\HH(Z)$, may be
thought of as a structural kind of resolution or ``chain complex'' by
adding suitable signs to the boundary operator.  $\HH(Z)$ may be
thought of as a structure on $Z$ or a structure measure.  These
structures may lead to new cohomological notions beyond $\PPP$ and
$\G$.

If the $\Shyper_i$'s have a tensor product it is natural to require:
\begin{center}
  \begin{tikzpicture}
    \diagram{d}{2.5em}{2.5em}{
      \hat{\delta}_{n - i + 1} \colon \otimes \{\Omega_{i - 1}
      (\{b_{i - 1}\})\} & \Omega_i(b_i)\\
    };

    \path[->] (d-1-1) edge (d-1-2);
  \end{tikzpicture}
\end{center}
relating our globalizers to factorization algebras and quantum field
theories (\citeasnoun{SO}).  The breakdown of a globalizer may lead to
disorganized states.  In cells in an organism this may be related to
cancerous states, see also \citeasnoun{HAM}.  Globalizers represent an
organizational way to get global properties and states of general
collections of objects.

In hyperstructures the release of properties, energy or creation of
states at one level may facilitate new releases or creations at the
higher levels.  This will continue until often an optimum is reached
at the top (global) level.  This is an analogy with the gluing of
states in Globs.  The mathematical details of this construction will
be followed up elsewhere.  Our main purpose here is to show how
hyperstructures may be used in obtaining global states, properties,
fields or objects from local ones.

\section{Bonds and bridges}
\label{sec:bridges}

Often in mathematics and science one would like to transfer a
structure in one context to a similar structure in another context.
Hyperstructures are useful tools in such situations transfering
structures via higher bonds.  This has been discussed in
\citeasnoun{SO}, but we will elaborate on the theme.

Putting a hyperstructure on a collection may be useful in the study
and use of the collection.  Sometimes one may transfer the structure
to another situation and use it there in order to get new results
otherwise difficult to detect or prove.

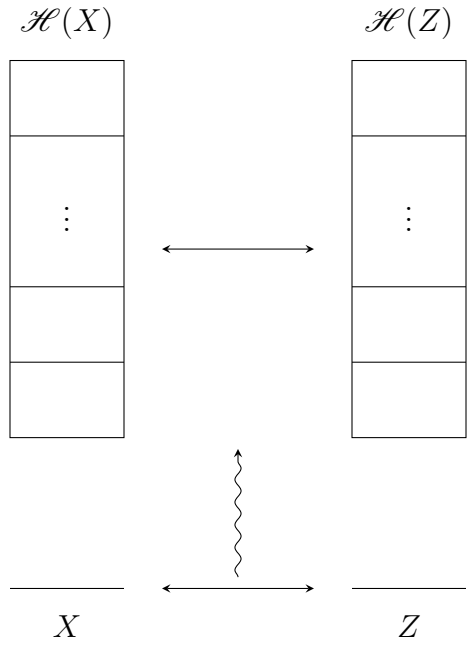
\begin{figure}[H]
  \centering
  \begin{tikzpicture}
    \draw (0,0) rectangle (1.5,5);
    \draw[<->] (2,2.5) -- (4,2.5);
    \draw (4.5,0) rectangle (6,5);

    \foreach \a in {1,2,4}{
      \draw (0,\a) -- (1.5,\a);
      \draw (4.5,\a)  -- (6,\a);
    }

    \draw (0,-2) -- (1.5,-2);
    \draw[<->] (2,-2) -- (4,-2);
    \draw (4.5,-2) -- (6,-2); 

    \draw[decorate, decoration={snake, amplitude=.4mm,
        segment length=3mm, post length=1mm},->] (3,-1.85) -- (3,-0.15);

    \node at (0.75,5.5){$\HH(X)$};
    \node at (5.25,5.5){$\HH(Z)$};
    \node at (0.75,-2.5){$X$};
    \node at (5.25,-2.5){$Z$};
    \node at (0.75,3){$\vdots$};
    \node at (5.25,3){$\vdots$};
  \end{tikzpicture}
  \caption{Structural transfer}
  \label{fig:structural_transfer}
\end{figure}

If we have a collection $X$ and a hyperstructure on it, $\HH(X)$, and
$X$ is ``related'' in some way to another collection $Z$, we may let
the relations propagate throught $\HH$ in order to form a new
hyperstructure $\HH(Z)$ as discussed in \citeasnoun{SO}, see Figure
\ref{fig:structural_transfer}.

\begin{figure}[H]
  \centering
  \begin{tikzpicture}
    \draw (0,0) rectangle (1.5,4);
    \draw (4.5,0) rectangle (6,4);

    \draw[<->] (2,2) -- node[above,midway]{bond} (4,2);

    \node at (3,4.5){$\HH$};
    \node at (0.75,2){$\HH_1$};
    \node at (5.25,2){$\HH_2$};

    \draw (-0.25,-0.25) rectangle (6.25,4.75);
  \end{tikzpicture}
  \caption{A hyperstructure bridge}
  \label{fig:bridge}
\end{figure}
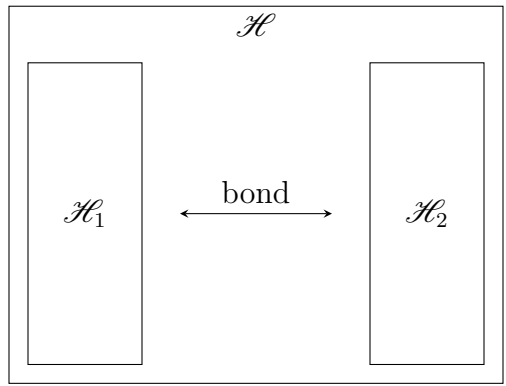

This induced transfer of structure is essentially a new kind of bond
in a hyperstructure embracing the two ($\HH(X)$ and $\HH(Z)$) as
indicated in Figure \ref{fig:bridge}.

\begin{figure}[H]
  \centering
  \subfigure{
    \centering
    \begin{tikzpicture}[font=\tiny]
      \begin{scope}
        \draw (0,0) rectangle (2.5,5);
        \foreach \y in {1,2,4}{
          \draw (0,\y) -- (2.5,\y);
        }

        \node at (1.25,-0.6){$\HH_{R_0}(Z)$};
        \node at (1.25,0.4){$\HH_{R_0}^0(Z) = Z$};
        \node at (1.25,1.4){$\HH_{R_0}^1(Z)$};
        \node at (1.25,3){$\vdots$};
        \node at (1.25,4.4){$\HH_{R_0}^n(Z)$};

        \draw[->] (3,0.4) -- node[midway,above]{$R_0$} (4.5,0.4);
      
        \foreach \y in {1.5,4.5}{
          \draw[dashed,->] (3,\y) -- (4.5,\y);
        }
      \end{scope}

      \begin{scope}[xshift=5cm]
        \draw (0,0) rectangle (2.5,5);
        \foreach \y in {1,2,4}{
          \draw (0,\y) -- (2.5,\y);
        }
        
        \node at (1.25,-0.6){$\HH(X_0)$};
        \node at (1.25,0.4){$\HH^0(X_0) = X_0$};
        \node at (1.25,1.4){$\HH^1(X_0) = X_1$};
        \node at (1.25,3){$\vdots$};
        \node at (1.25,4.4){$\HH^n(X_0) = X_n$};
        
        \draw[->] (3,0.4) -- node[midway,above]{$I_0$} (4.5,0.4);
        
        \foreach \y in {1.5,4.5}{
          \draw[dashed,->] (3,\y) -- (4.5,\y);
        }
      \end{scope}

      \begin{scope}[xshift=10cm]
        \draw (0,0) rectangle (2.5,5);
        \foreach \y in {1,2,4}{
          \draw (0,\y) -- (2.5,\y);
        }
      
        \node at (1.25,-0.6){$\HH_{I_0}(Z)$};
        \node at (1.25,0.4){$\HH_{I_0}^0(Z) = Z$};
        \node at (1.25,1.4){$\HH_{I_0}^1(Z)$};
        \node at (1.25,3){$\vdots$};
        \node at (1.25,4.4){$\HH_{I_0}^n(Z)$};
      \end{scope}
    \end{tikzpicture}
  } \subfigure{
    \begin{tikzpicture}[scale=0.7,font=\small]
      \node at (-1,0){};

      \begin{scope}
        \draw (0,0) rectangle (5,9);
        \foreach \y in {3,6}{
          \draw (0,\y) -- (5,\y);
        }

        \begin{scope}[xshift=2.5cm,yshift=7.5cm,scale=0.55]
          \draw[very thick] (0,0) circle(2.5cm);
        
          \begin{scope}[scale=0.3,xshift=-4.25cm,yshift=2.5cm]
            \draw[very thick] (0,0) circle(2.5cm);
            \draw[very thick, blue] (-0.25,-1.6) circle(0.5cm);
            \draw[very thick, red] (-1.1,1.2) circle(0.5cm);
            \draw[very thick, green] (1.4,0.8) circle(0.5cm);
          \end{scope}
        
          \begin{scope}[scale=0.3,xshift=4.5cm,yshift=1.2cm]
            \draw[very thick] (0,0) circle(2.5cm);
            \draw[very thick, blue] (-0.25,-1.6) circle(0.5cm);
            \draw[very thick, red] (-1.1,1.2) circle(0.5cm);
            \draw[very thick, green] (1.4,0.8) circle(0.5cm);
          \end{scope}
          
          \begin{scope}[scale=0.3,xshift=-0.5cm,yshift=-4cm]
            \draw[very thick] (0,0) circle(2.5cm);
            \draw[very thick, blue] (-0.25,-1.6) circle(0.5cm);
            \draw[very thick, red] (-1.1,1.2) circle(0.5cm);
            \draw[very thick, green] (1.4,0.8) circle(0.5cm);
          \end{scope}
        \end{scope}
        
        \begin{scope}[xshift=2.5cm,yshift=4.5cm,scale=0.65]
          \begin{scope}[scale=0.3,xshift=-4.25cm,yshift=2.5cm]
            \draw[very thick] (0,0) circle(2.5cm);
            \draw[very thick, blue] (-0.25,-1.6) circle(0.5cm);
            \draw[very thick, red] (-1.1,1.2) circle(0.5cm);
            \draw[very thick, green] (1.4,0.8) circle(0.5cm);
          \end{scope}
        
          \begin{scope}[scale=0.3,xshift=4.5cm,yshift=1.2cm]
            \draw[very thick] (0,0) circle(2.5cm);
            \draw[very thick, blue] (-0.25,-1.6) circle(0.5cm);
            \draw[very thick, red] (-1.1,1.2) circle(0.5cm);
            \draw[very thick, green] (1.4,0.8) circle(0.5cm);
          \end{scope}
        
          \begin{scope}[scale=0.3,xshift=-0.5cm,yshift=-4cm]
            \draw[very thick] (0,0) circle(2.5cm);
            \draw[very thick, blue] (-0.25,-1.6) circle(0.5cm);
            \draw[very thick, red] (-1.1,1.2) circle(0.5cm);
            \draw[very thick, green] (1.4,0.8) circle(0.5cm);
          \end{scope}
        \end{scope}
        
        \foreach \x in {1.25,2.5,3.75}{
          \draw[thick,red] (\x,2.5) circle(0.3cm);
          \draw[thick,green] (\x,1.5) circle(0.3cm);
          \draw[thick,blue] (\x,0.5) circle(0.3cm);
        }

        \node at (2.5,9.8){State hyperstructure:};
        \node at (2.5,-0.6){$Z = \text{$9$ particle states}$};
        
        \foreach \y in {4.5,7.5}{
          \draw[line width=1.5,->] (6,\y) -- (9,\y);
        }

        \draw[line width=1.5,->] (6,1.5) -- node[midway,below]{$R_0$}
        (9,1.5);
      \end{scope}
      
      \begin{scope}[xshift=10cm]
        \draw (0,0) rectangle (5,9);
        \foreach \y in {3,6}{
          \draw (0,\y) -- (5,\y);
        }
        
        \begin{scope}[knot/.style={thin knot},scale=0.25]
          \begin{scope}[xshift=5cm,yshift=20cm]
            \setbrstep{.1}
            \brunnian{2.5}{3}
          \end{scope}
          
          \begin{scope}[xshift=15cm,yshift=20cm]
            \setbrstep{.1}
            \brunnian{2.5}{3}
          \end{scope}
        
          \begin{scope}[xshift=10cm,yshift=16cm]
            \setbrstep{.1}
            \brunnian{2.5}{3}
          \end{scope}
        \end{scope}
        
        \foreach \x in {1.25,2.5,3.75}{
          \draw[thick,red] (\x,2.5) circle(0.3cm);
          \draw[thick,green] (\x,1.5) circle(0.3cm);
          \draw[thick,blue] (\x,0.5) circle(0.3cm);
        }
        
        \node at (2.5,9.8){Link hyperstructure:};
        \node at (2.5,-0.6){$X_0 = \text{$9$ rings}$};
      \end{scope}
    \end{tikzpicture}
    \raisebox{4.8cm}[0pt][0pt]{
      \hspace*{-3.175cm}
      \includegraphics[scale=0.09]{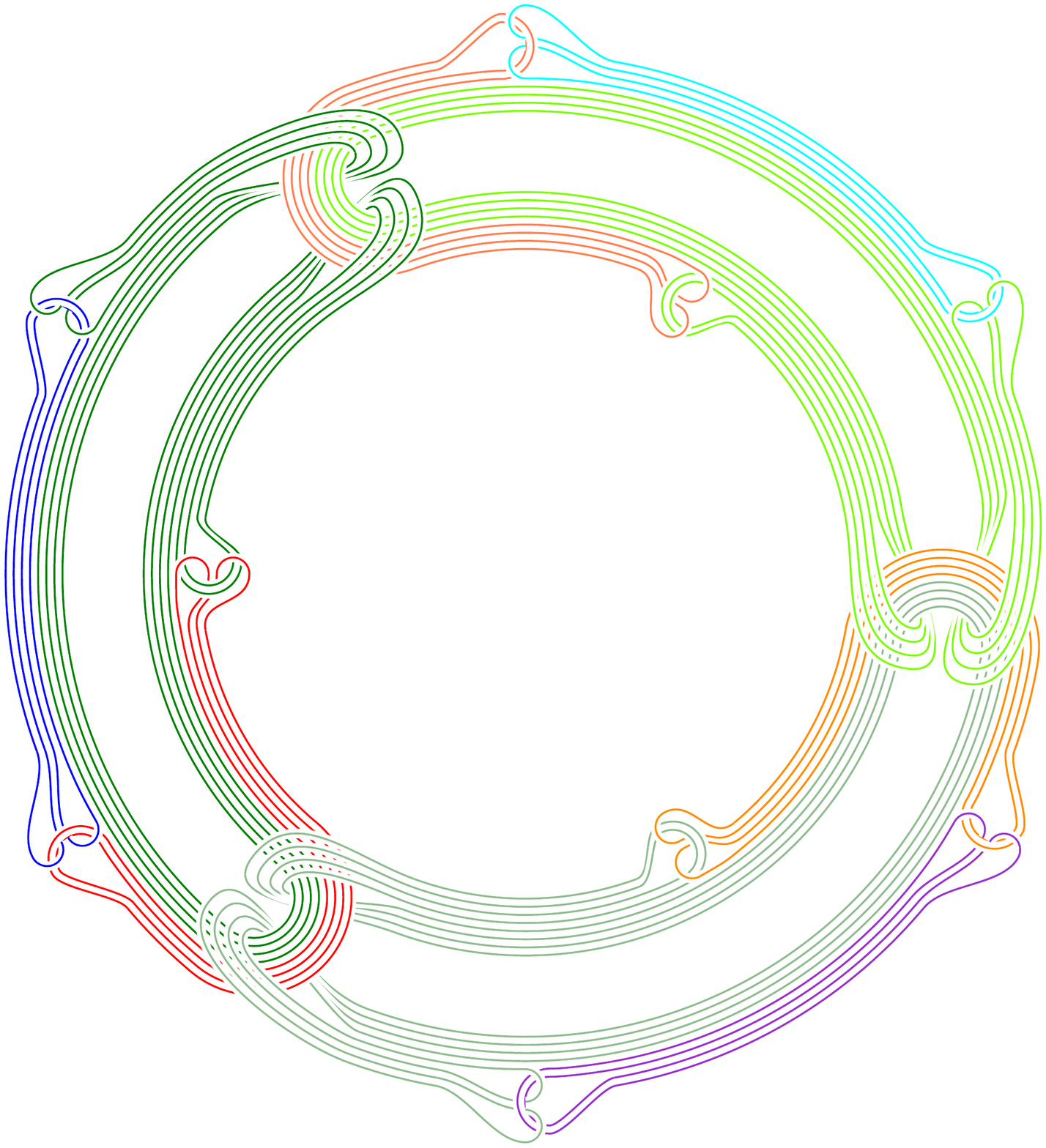}
    }
  } \subfigure{
  }
  \caption{Transfer of structures via hyperstructures}
  \label{fig:induced_transfer}
\end{figure}

Such a bond will act as a bridge between the two structures
transferring at the various levels.

An example of such transfer of structure is the following discussed in
\citeasnoun{SO}:

A finite set of rings in three dimensional space may be given an
$\HH$-structure with a levelwise Brunnian property: the rings at any
level are linked such that if one is removed, the rest become
unlinked.  Then we may consider a (quantum) particle system and by
representing particles by rings we may give the particle system a
similar $\HH$-structure with the Brunnian property: if at any level a
certain number of particle collections interact in such a way that if
one collection is removed, the remaining collections do not interact.
The representation
\begin{center}
  \begin{tikzpicture}
    \diagram{d}{2.5em}{2.5em}{
      \mathrm{ring} & \mathrm{particle}\\
    };

    \path[->] (d-1-1) edge (d-1-2);
  \end{tikzpicture}
\end{center}
induces the bridge of hyperstructures, see Figure
\ref{fig:induced_transfer}.

Bridges have been discussed in many interesting ways by
\citeasnoun{Caramello}.  The main idea being transfer of information
between two different sites via categories or toposes of sheaves as
follows: given two sites
\begin{equation*}
  (\C_1,J_1) \quad \text{and} \quad (\C_2,J_2)
\end{equation*}
We may have the situation
\begin{center}
  \begin{tikzpicture}
    \diagram{d}{2.5em}{2.5em}{
      (\C_1,J_1) & \Sh(\C_1,J_1) \simeq \Sh(\C_2,J_2) & (\C_2,J_2)\\
    };

    \path[->]
      (d-1-1) edge (d-1-2)
      (d-1-2) edge (d-1-3);
  \end{tikzpicture}
\end{center}
where $\simeq$ means (Morita-type) equivalence.  If such an
equivalence exists it represents a bridge for transfer of information
between $(\C_1,J_1)$ and $(\C_2,J_2)$.  In our terminology we think of
the equivalence as a bond in the hyperstructure or the category of
Grothendieck toposes.

In the hyperstructure setting this corresponds to given two
$\HH$-sites $(\HH_1,J_1)$ and $(\HH_2,J_2)$ with a bridge or bond of
the globalizer hyperstructures:
\begin{center}
  \begin{tikzpicture}
    \diagram{d}{2.5em}{2.5em}{
      (\HH_1,J_1) & \Glob(\HH_1,J_1) \simeq \Glob(\HH_2,J_2) &
      (\HH_2,J_2)\\
    };

    \path[->]
      (d-1-1) edge (d-1-2)
      (d-1-2) edge (d-1-3);
  \end{tikzpicture}
\end{center}
where $\simeq$ means a suitable equivalence or bond.

In general when
\begin{center}
  \begin{tikzpicture}
    \diagram{d}{0.5em}{3em}{
      (\HH_1,J_1) &\\
      & (\HH,J)\\
      (\HH_2,J_2) &\\
    };

    \path[right hook->]
      (d-1-1.east) edge (d-2-2.north west)
      (d-3-1.east) edge (d-2-2.south west);
  \end{tikzpicture}
\end{center}
one may transfer structures via the top level in $\HH$ and similarly
for states via
\begin{center}
  \begin{tikzpicture}
    \diagram{d}{0.5em}{3em}{
      \Glob(\HH_1,J_1) &\\
      & \Glob(\HH,J).\\
      \Glob(\HH_2,J_2) &\\
    };

    \path[right hook->]
      (d-1-1.east) edge (d-2-2.north west)
      (d-3-1.east) edge (d-2-2.south west);
  \end{tikzpicture}
\end{center}

Another possible type of bridge using hyperstructures would be to
introduce a ``Brave New $\HH$-algebra'' where equalities (and hence
equations) are being replaced by bonds in hyperstructures.  This
extends homotopical algebra (homotopy bonds) and gives many new
perspectives to be explored.

\section{Organizing collections}
\label{sec:org_collections}

When working with collections of objects it is often desirable to make
new collections from old ones.  One may break up collections to new
ones, break up objects to more objects (fission) or amalgamate objects
and collections (fusion).  As we have pointed out organizing
collections is useful in detecting and creating properties and
structure, controlling actions and making constructions, classifying
and achieving goals.  Our thesis is \emph{that hyperstructured
  architectures are very useful in doing this}.

In this setting fusion and fission are dual.  In a fusion process one
takes the given collection as the lowest level, and in a fission
process we let the collection represent the top level.  This is
analogous to preparing a collection for analysis or synthesis.  Our
collections could consist of any types of objects, be sets, spaces
(topological and metric spaces, manifolds, simplicial spaces,
varieties,$\ldots$), groupoids, even hyperstructures themselves.  To
illustrate our ideas we will here basically consider fusion-type
processes.

\begin{figure}[H]
  \centering
  \includegraphics{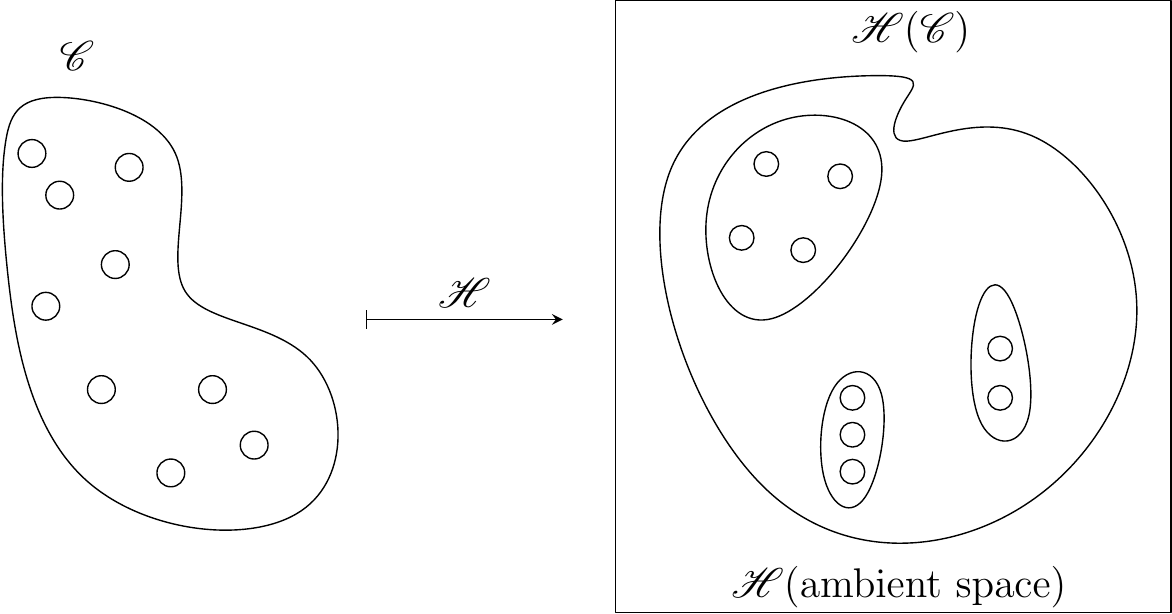}
  \caption{Fusion of a collection}
  \label{fig:fusion}
\end{figure}

\begin{figure}[H]
  \centering
  \includegraphics{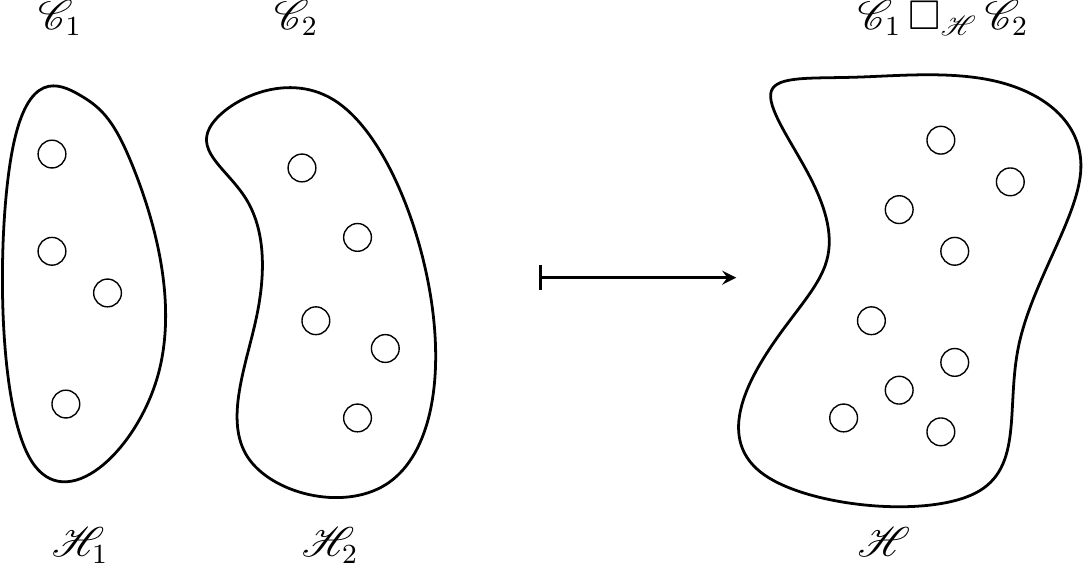}
  \caption{Fusing collections}
  \label{fig:fusion_collections}
\end{figure}

Putting a hyperstructure on a collection in such a way that the
collection is represented as the lowest level of the hyperstructure is
essentially a fusion process in the sense that objects are levelwise
being bound.  Then local states and properties may propagate through
the hyperstructure by globalizers to the top level.  This represents a
fusion of existing objects into collections of collections etc.\
determinded by the bonds.  Such a fusion process may be induced by a
hyperstructure on the collection or a hyperstructure on the ambient
space of the collection which again will induce a hyperstructure on
the collection.  This may be illustrated as in Figure \ref{fig:fusion}.
The next basic step is the fusion of a collection of collections into
a collection (Figure \ref{fig:fusion_collections}.  Let us just
consider two collections $\C_1$ and $\C_2$, the process defines a
``fusion product'':
\begin{center}
  \begin{tikzpicture}
    \diagram{d}{2.5em}{2.5em}{
      (\C_1,\C_2) & \C_1 \, \fusionscript \, \C_2\\
    };

    \path[|->] (d-1-1) edge (d-1-2);
  \end{tikzpicture}
\end{center}

The new thing to study is that $\fusionscript$ depends on a
hyperstructure ($\HH$-structure) on $\C_1$ and $\C_2$ in order to
obtain the fusion.  They may possibly be different and $\fusionscript$
stands as a generic notion.  States and properties may then be built
in such that $\C_1 \fusionscript \C_2$ is a collection with a derived
global state or property
\begin{equation*}
  \omega (\C_1 \fusionscript \C_2)
\end{equation*}
determined from $\omega(\C_1)$ and $\omega(\C_2)$ via globalizers in
$\HH$.  Such a process may reflect a ``release'' of properties of the
bonds --- like for example a change from high to low energy levels.

The second type of fusion to be considered is when objects fuse into
new types of objects, hence we form a new collection $\C$ from $\C'$:
\begin{center}
  \begin{tikzpicture}
    \diagram{d}{2.5em}{2.5em}{
      \C & \C'\\
    };

    \path[|->] (d-1-1) edge (d-1-2);
  \end{tikzpicture}
\end{center}
\begin{figure}[H]
  \centering
  \includegraphics{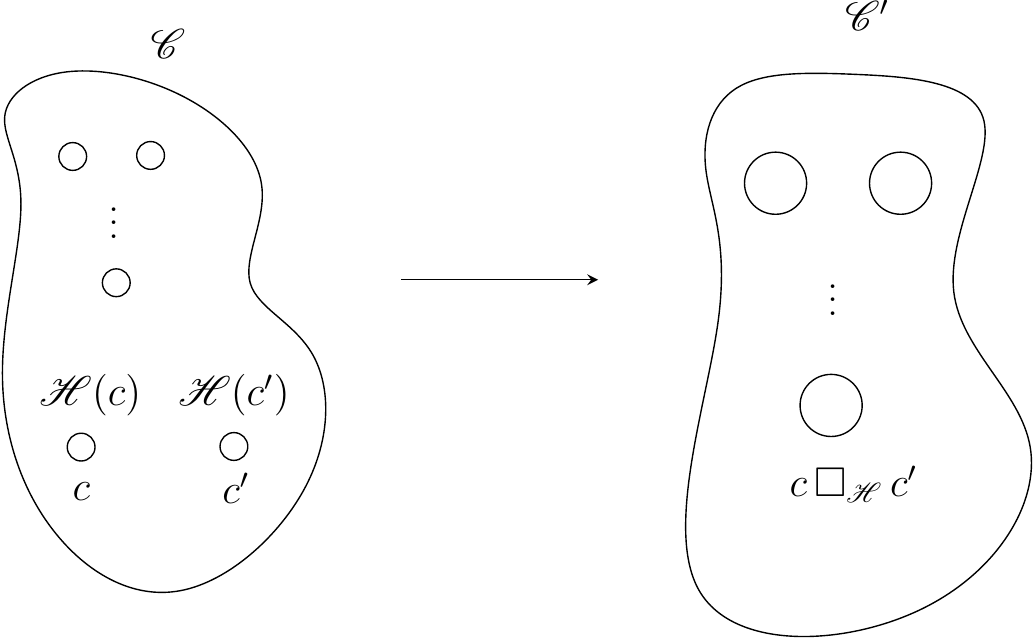}
  \caption{Object fusion in a collection by object hyperstructures}
  \label{fig:object_fusion}
\end{figure}
When such a process is desired one may assign a hyperstructure
$\HH(c)$ to each object in the collection, $c\in \C$ (Figure
\ref{fig:object_fusion}).  $\HH(c)$ may be induced from a
hyperstructure $\HH_A$ of an ambient space $A$.  There may also be an
$\HH(\C)$ on the total collection and they may possibly all interact
in order to create globalizers with desired properties (for example
like release of energy as in nuclear fusion).  This applies not only
to fusion of pairs, but to arbitrary collections (finite or infinite).

Fusion processes of this type give rise to many interesting aspects
like introducing ``energy'' thresholds for fusion at each level and
studying how the fusion process propagates through the $\HH$-organized
collection.

Our main point of view is that a hyperstructure on objects,
collections or ambient space may facilitate fusion (and similar)
processes, being an application of (see \citeasnoun{SO}):

\begin{thp}
  \emph{In order to study or use a collection of objects it is useful
    to put on a suitable Hyperstructure to reach a goal.
    Hyperstructures are tools for thought and tools for organizing
    complex information.}
\end{thp}

\section{Types, bonds and hyperstructures}
\label{sec:types}

Type theory is a formal logical system intended to provide a framework
for mathematical thinking and constructions.  It has also turned out
to be useful as a framework for high-level programming languages.  In
traditional type theory one may think of types as sets, but in a new
approach V.\ Voevodsky and co-workers have introduced Homotopy Type
Theory where one basically thinks of types as spaces and homotopy
types, Voevodsky \& co \cite{HTT}.  Furthermore, to each type
one may associate an $\infty$-groupoid.  Voevodsky advocates building
mathematics on ($\infty$-)groupoids instead of sets.

As pointed out in \cite{Cog,NS,SO} we think of hyperstructures as in
many ways reflecting evolutionary processes and the way we think and
make things.  Shortly they represent: Tools for Thought.  This becomes
quite apparant in higher order thinking when we have to consider
several levels at a time.  It may be that the human brain is basically
designed to consider three levels at a time: where we are, the level
above and the level below, see for example \citeasnoun{Koestler}.
Beyond that we need a framework in which to reason, since
comprehending higher order structures and hierarchies is vital in
human thinking.  We clearly see this in axiomatic systems (ZFC) for
set theory.

In a lecture at IAS in Princeton March 26, 2014, V.\ Voevodsky
suggested that a system adequate both for human reasoning and computer
verifications should contain:
\begin{itemize}
\item[i)] a formal deduction system
\item[ii)] providing meaning comprehensible to humans
\item[iii)] a way for humans to encode mathematical ideas (like
  hierarchies).
\end{itemize}

We suggest that hyperstructures do this and may be extended to a new
type theory.  It is desirable that deductions are directly built into
the system.  This is the case for hyperstructures as we will show.

The simple idea is: Think of bonds as types.  We have higher bonds
(bonds of bonds $\ldots$) giving higher types.  So higher order
structures and hierarchies are built in.  The semantics (meaning
aspect) is built in via the $\Omega$ assignments (``states or
properties'').

The deducational part goes as follows: a bond $b_1$ at one level binds
bonds at a lower level: $\{a_o^i\}$, $\{b_0^j\}$.  

Assume that for every bond $b_1$ and a given collection $\{a_o^i\}$
one can decide whether
\begin{equation*}
  b_1(\{a_o^i\},\{b_o^j\})
\end{equation*}
holds for the collection $\{b_o^j\}$.  If so, we say that $\{b_o^j\}$
can be deduced from $\{a_o^i\}$.

One may even allow an intermediate ``proof'' --- family $\{c_0^k\}$,
such that
\begin{equation*}
  b_1(\{a_0^i\},\{c_0^k\},\{b_0^j\}).
\end{equation*}
Geometrically this may be thought of as a filling-in condition a la
Kan-conditions in homotopy theory.  Oriented bonds may also be
considered.

If the bonds are represented as spaces one may intuitively think of a
proof as filling in the missing pieces of a decomposition (for example
in a triangulation of a space).  This will allow for a geometric way
of thinking of deductions.

Universes of bonds of order $n$ will be bonds of order $(n + 1)$.
Since each bond governs a sub-hyperstructure one will have the
association of a type = bond = $A$, to a hyperstructure $\HH_A$,
similar to the $\infty$-groupoid association in homotopy type theory.

In the deduction both syntax and semantics play a role.  Comprehension
of (global) meaning comes from the existence of globalizers (global
type sheaves).  This is compatible with our point of view of
hyperstructures as also describing evolutionary structures.  For
example in a biological system --- an organism  --- a globalizer would
assign a global biological state --- useful in taking actions.

In hyperstructure types one may have a variety of ``identity types'',
\begin{equation*}
  b_1(b_0^{i_1}), \quad b_1(b_0^{i_1},b_0^{i_2}), \quad \ldots \quad, \quad
  b_1(b_0^{i_1},\ldots,b_0^{i_n})
\end{equation*}
--- if needed, motivated by geometric cobordisms.

Hyperstructures as introduced here are based on sets or collections.
They exist in the universe of sets with structure.  However, in the
universe of homotopy type theory one may study $\infty$-groupoids with
hyperstructures added on as additional structure.  Or if a type theory
of hyperstructures is taken as a basis for our thinking we would
consider hyperstructures with other imposed hyperstructures.  In
general, whichever basic building blocks we choose for our formal
thinking system, we should also in addition consider them equipped
with hyperstructures.

In order to describe, encode and manipulate (mathematical) ideas in a
formal system symbols are needed.

Meaningful symbol manipulations result in collections of symbols
organized into a hyperstructure which may improve the encoding of many
mathematical ideas.  This allows for a symbol representation by higher
dimensional shapes, not just one dimensional shapes and one
dimensional concatenation.  Symbols may be higher dimensional shapes
and concatenation being replaced by higher order gluing in the form of
bonds of symbols being bonds of bonds $\ldots$ Symbol manipulation is
then controlled by higher order bonds and rules for bond composition.
We will return to these issues in future papers.

\section{Conclusion}
\label{sec:conclusion}

Using hyperstructures gives rise to a variety of higher order
architectures on general collections of objects.  We have argued that
putting a hyperstructure on a collection or situation may be very
useful, and new concepts have been introduced.  Higher order gluing,
bridges, fusion of collections and local to global situations have
been discussed.

We think that the higher order architectures based on hyperstructures
may in the future turn out to be very useful in many areas of science:
mathematics, physics, chemistry, biology and neuroscience, computer
science, economics and social science, engineering, systems science
and architecture itself, as a matter of fact in all areas of human
thought.

\section*{Acknowledgements}
I would like to thank Marius Thaule for valuable technical assistance
with the manuscript, and Andrew Stacey for the graphical production of
Figure \ref{fig:links} and parts of Figure \ref{fig:induced_transfer}.
The rest of the figures were made by Marius Thaule.  Furthermore, I
would also like to thank The Institute for Advanced Study, Princeton,
NJ, USA for their kind hospitality during my stay there in 2013 when
parts of this work were done.



\end{document}